\documentclass[conference]{IEEEtran}

\usepackage[T1]{fontenc}
\usepackage{subcaption}
\usepackage{hyperref}
\usepackage{amsmath}
\usepackage{amssymb}
\usepackage{accents}
\usepackage{mathtools}
\newcommand{\ubar}[1]{\underaccent{\bar}{#1}}
\usepackage{cleveref}
\usepackage{graphicx}
\usepackage{lipsum}
\usepackage{booktabs}
\usepackage{array}
\usepackage{comment}
\usepackage{cite}
\usepackage{multirow}
\usepackage[acronym]{glossaries}

\makeglossaries

\newacronym{pv}{PV}{Photovoltaic}
\newacronym{pdf}{PDF}{Probability Density Function}
\newacronym{pmf}{PMF}{Probability Mass Function}
\newacronym{cdf}{CDF}{Cumulative Distribution Function}
\newacronym{ds}{DiS}{Dispatch Schedule}
\newacronym{bess}{BESS}{Battery Energy Storage System}
\newacronym{ev}{EV}{Electric Vehicle}

\begin{document}

\title{Probabilistic Day-Ahead Battery Scheduling based on Mixed Random Variables for Enhanced Grid Operation}
\author{Janik Pinter\IEEEauthorrefmark{1}, Frederik Zahn\IEEEauthorrefmark{1}, Maximilian Beichter\IEEEauthorrefmark{1}, Ralf Mikut\IEEEauthorrefmark{1}‚ and Veit Hagenmeyer\IEEEauthorrefmark{1}
\\Institute for Automation and Applied Informatics‚
\\Karlsruhe Institute of Technology (KIT), Germany
\\ \IEEEauthorrefmark{1}{\tt\small \{firstname.lastname\}@kit.edu}
}

\maketitle

\begin{abstract}
The increasing penetration of renewable energy sources introduces significant challenges to power grid stability, primarily due to their inherent variability. A new opportunity for grid operation is the smart integration of electricity production combined with battery storages in residential buildings. This study explores how residential battery systems can aid in stabilizing the power grid by flexibly managing deviations from forecasted residential power consumption and PV generation. The key contribution of this work is the development of an analytical approach that enables the asymmetric allocation of quantified power uncertainties between a residential battery system and the power grid, introducing a new degree of freedom into the scheduling problem. This is accomplished by employing mixed random variables---characterized by both continuous and discrete events---to model battery and grid power uncertainties. These variables are embedded into a continuous stochastic optimization framework, which computes probabilistic schedules for battery operation and power exchange with the grid. Test cases demonstrate that the proposed framework can be used effectively to reduce and quantify grid uncertainties while minimizing electricity costs. It is also shown that residential battery systems can be actively used to provide flexibility during critical periods of grid operation. Overall, this framework empowers prosumers to take an active role in grid stabilization, contributing to a more resilient and adaptive energy system.
\end{abstract}


\section{Introduction}
\label{sec:introduction}
The global transition toward renewable energy generation is transforming electrical power systems worldwide. While renewable energy sources offer significant environmental and economic benefits, they also introduce new challenges due to their intrinsic volatility. Historically, uncertainties in power systems primarily occurred on the demand side and were effectively manageable due to the controllability of conventional power plants, i.e., the generation was able to match the demand. As the penetration of renewable energy sources grows, the controllability on the generation side decreases, leading to a significant increase in overall uncertainties within the power grid. To enable a secure and reliable grid operation in the future, various methods for providing unidirectional flexibility and managing uncertainties at different power and voltage levels are currently being investigated, such as, at
\begin{itemize}
    \item Grid Level: Large-scale balancing mechanisms and grid-level storage systems are investigated to ensure system-wide stability \cite{girigoudar2022chance}.
    \item Subgrid Level: Microgrids coordinate local renewable generation, energy storages, and demand-side management to handle uncertainties and reduce dependence on the main grid \cite{konneh2022application}. 
    \item Single Unit: Individual industrial or residential buildings can decrease their power consumption uncertainties by optimized scheduling of energy consumption and \gls{bess} operation \cite{sossan2016achieving}.
\end{itemize}
The increasing deployment of residential \gls{pv} systems in conjunction with \glspl{bess} demands further investigation on how to best integrate them into power system management. The major sources of uncertainty on a household level are deviations from forecasted consumption profiles and forecasted \gls{pv} production. \glspl{bess} offer great potential to balance these deviations, as opposed to injecting them into the grid. To achieve this, the operation of \glspl{bess} need to take into account the uncertainties in a systematic way. Therefore, the present work focuses on the development of a scheduling approach for a single house, which allows to reduce and quantify grid power injection uncertainties. We present a novel formulation to calculate an optimized day-ahead schedule for \gls{bess} operation in residential buildings featuring a \gls{pv} system.


\subsection{Related Works}
Several studies have proposed methods for optimizing day-ahead \gls{bess} scheduling in residential buildings, which differ in two distinct aspects: First, how uncertainties — whether in generation, demand, or electricity prices — are embedded into the optimization process. Second, the specific objective the optimization seeks to achieve.

Uncertainties can be embedded into an optimization framework in multiple ways. One common approach is Robust Optimization (RO), which focuses on optimizing system performance while ensuring feasibility for worst-case realizations. Typically, uncertain parameters are assumed to lie within fixed bounds \cite{bahramara2021robust, haider2023robust, cao2020optimal, ebrahimi2019adaptive}. The optimization problem is solved using their worst-case values, i.e., most often the boundary values of the assumed intervals \cite{xie2024distributionally}. For example, \cite{bahramara2021robust} optimizes the day-ahead scheduling of a smart home equipped with \gls{pv} and \gls{bess}, treating energy prices as uncertain within predefined bounds. 

Another popular approach for incorporating uncertainties into an optimization framework are scenario-generation techniques \cite{rezaeimozafar2024residential, wang2020scenario, tavakoli2023development}. Methods such as Sample-Average-Approximations, Markov Chains, and Monte-Carlo approaches create an ensemble of deterministic scenarios, each representing a possible realization of uncertain parameters. These methods transform the stochastic formulation into a more elaborate deterministic one, which allows the use of well-established, computationally efficient optimization techniques.
For instance, \cite{rezaeimozafar2024residential} minimizes residential electricity costs by combining a neural network with an error analysis technique to represent uncertain PV behavior in the form of a scenario ensemble.

However, a critical limitation of both RO and scenario-generation techniques is their inability to provide quantified uncertainty propagation throughout the optimization process. For our work, the focus is not just on optimal \gls{bess} scheduling but also on the power exchange with the grid. Therefore, we require a method that quantifies uncertainties in a coherent and probabilistic manner throughout the optimization process. Neither RO nor scenario-generation techniques meet this requirement, which is why we propose a different approach.

In our work, Stochastic Programming (SP) is employed to represent uncertainties in the optimization framework. SP allows for quantified uncertainty propagation by representing uncertain parameters as random variables with known probability distributions. These distributions can be incorporated into an optimization framework through chance constraints \cite{appino2018use, werling2022towards, singh2021lagrangian}, or by minimizing their expected values \cite{su2021optimal, scarabaggio2021distributed}. 

One limitation of SP, however, is that it assumes the probability distributions of random variables are fully known, which does not always reflect reality. This issue is addressed in Distributionally Robust Optimization (DRO) \cite{xie2024distributionally}. Instead of relying on a single probability distribution, DRO considers a set of potential probability distributions and optimizes for the worst-case distribution within the set \cite{shi2022day, parvar2022optimal, xu2020scheduling}. However, in our case, the determination of the worst-case distribution is not straightforward, because the expected values of the random variables are not fixed but depend on the decision variables, determined during optimization. Thus, while DRO presents an appealing theoretical framework, its practical application is limited for the given case.


In general, cutting costs is the main incentive to optimize \gls{bess} operation. This is most commonly achieved by performing peak shaving \cite{su2022optimization}, load shifting \cite{sharda2021demand}, price arbitrage operations \cite{su2021optimal} or maximizing self-sufficiency \cite{ciocia2021self}.
Our study investigates the extension of those conventional usages.
The fundamental idea, first presented in \cite{sossan2016achieving}, is that the homeowner communicates its forecasted grid power exchange to the grid operator one day in advance. In the following, this forecasted power exchange is called \gls{ds}. The \gls{bess} is used to minimize the deviations from the \gls{ds}, without neglecting the \gls{bess}'s original purpose of minimizing electricity costs. By doing so, the overall uncertainty of the power grid can be reduced, which leads to a decrease in balancing power requirements and system costs \cite{mlilo2021impact}.
The framework from \cite{sossan2016achieving} is divided into two stages: the day-ahead stage, where the \gls{ds} is computed, and the real-time stage, where Model Predictive Control is applied to minimize \gls{ds} deviations.

The authors of \cite{appino2018use} enhance the proposed framework by incorporating probabilistic power forecasts. They assume that the \gls{bess} is able to fully compensate for all deviations from forecasted consumption and PV production profiles. As a result, the \gls{ds} remains deterministic, whereas the \gls{bess} schedule is probabilistic. To ensure feasibility, the \gls{bess} constraints, i.e., its power and capacity limits, are enforced through chance-constraints. These constraints ensure compliance with a given probability, known as security level. This approach focuses on a reliable \gls{ds} by shifting all uncertainties toward the \gls{bess}. 
However, this limits the ability to minimize electricity costs by exploiting the residential \gls{bess}.
Additionally, selecting appropriate security levels and interpreting the individual chance-constraints are challenging themselves \cite{singh2023chance}, and real-world limitations, like \gls{bess} response times or ramping constraints, prevent the full compensation of \gls{ds} deviations by the \gls{bess} on smaller time scales \cite{beichter2023towards}.





The authors of \cite{su2021optimal} investigate a similar probabilistic setting, but focus on minimizing residential electricity costs by formulating an exact expectation-based stochastic optimization approach. Convexity of the nonlinear optimization problem is achieved and proven. The uncertainty is handled by assuming that the grid fully compensates for all power uncertainties. Therefore, in contrast to \cite{appino2018use}, the approach results in a probabilistic \gls{ds} and a deterministic \gls{bess} schedule. In the latter case, the \gls{ds} is more of a by-product instead of a fundamental objective. 

To summarize, SP approaches for residential day-ahead \gls{bess} scheduling tend to either assign all uncertainties to the grid \cite{su2021optimal}, or to the \gls{bess} \cite{appino2018use}, sidestepping the complexity of allocating uncertainties between the two systems. Approaches where the \gls{bess} absorbs all uncertainties support the power grid but face challenges in interpretability and electricity cost minimization.
Conversely, assigning all uncertainties to the grid allows for electricity cost minimization, but disregards the inclusion of residential \glspl{bess} to provide flexibility to support the power grid.   
Thus, there are no scheduling methods that can balance both systems’ roles and constraints. This gap is filled in the present work.



\subsection{Contribution}

The main idea of the present work is a formulation of the \gls{bess} scheduling problem, which allows to share quantified uncertainties between the grid and the \gls{bess}. 
The challenge lies in how to represent and manage these uncertainties effectively and how to incorporate them into a stochastic optimization framework. 
In the present work, we do not schedule set points for the \gls{bess}'s charging, but we assign intervals of minimum and maximum charging power. This formulation allows to absorb a limited amount of deviations in the \gls{bess} while respecting its physical limits and quantifying the uncertainty injected into the grid.

The key contributions of this paper include the following:
\begin{enumerate}
    \item Formulation of a novel approach to allow an asymmetrical allocation of quantified power uncertainties between a \gls{bess} and the power grid using mixed random variables.
    \item Derivation of an expectation-based nonlinear stochastic optimization problem to compute probabilistic schedules for \gls{bess} operation and grid exchange.
    \item Demonstration of the proposed scheduling algorithms for three different scenarios based on real-world data. 
\end{enumerate}

The paper is organized as follows: 
\Cref{sec:model-derivation} introduces the problem and contains the main contribution of this paper. The core idea on how to tackle uncertainties is presented and a mathematical power exchange model is formulated. In \Cref{sec:optimization-problem}, the derived model is embedded into an optimization framework in a generic form and is applied and analyzed in \Cref{sec:model-application}. \Cref{sec:conclusion} concludes this paper.

\section{Model Derivation}
\label{sec:model-derivation}

We investigate the combination of a single residential house, equipped with \gls{pv} panels and a \gls{bess}. We primarily focus on stationary \glspl{bess} but want to mention that the usage of \glspl{ev} with bidirectional charging capabilities essentially results in the same model, as long as the \glspl{ev} remain connected to the charging point. For \glspl{ev} that only allow unidirectional charging, the presented approach would shift more toward demand-side management. 

\begin{figure}[h]
    \centering
    \includegraphics[width=0.6\linewidth]{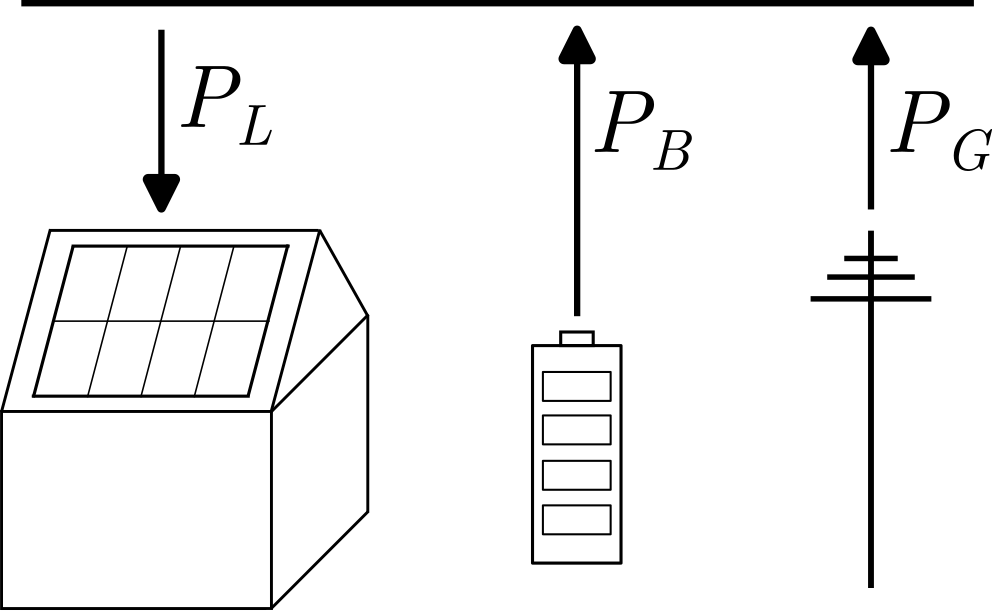}
    \caption{General setting. Residential prosumption, \gls{bess} power, and grid power are in balance: $P_L=P_B+P_G.$}
    \label{fig:problem-setting}
\end{figure}

In this study, the inflexible power demand (also referred to as load) is combined with the inflexible \gls{pv} power generation in the variable $P_L$. It represents both power production and consumption and is therefore referred to as prosumption. The basic structure is sketched in \autoref{fig:problem-setting}. $P_B$ corresponds to the controllable \gls{bess} power and $P_G$ denotes the power provided by the grid. In order to fully meet the residential power prosumption and under the assumption of a lossless power connection, the power exchange can be described via 
\begin{equation}
    \label{eq:power-balance}
    P_L = P_B + P_G,
\end{equation}
with positive power flows defined by the arrows in \autoref{fig:problem-setting}.

\subsection{Uncertainty Modeling}
\label{sec:forecasts}

In this work, the prosumption is modeled as a real-valued random variable $P_L\colon\Omega \rightarrow \mathbb{R}$ defined on a probability space $(\Omega, \mathcal{F}, \mathbb{P})$, where $\Omega$ is the set of possible outcomes, $\mathcal{F}$ is a sigma-algebra of measurable events, and $\mathbb{P}$ is a probability measure. $P_L$ represents the average power exchange over a specific time period\footnote{For the sake of simplicity, the time index is omitted in this subsection.} and is assumed to be absolutely continuous, implying that a \gls{pdf} $f_{P_L}$ exists.

We split the prosumption
\begin{equation}
    \label{eq:prosumption}
    P_L = \hat{p}_L + \Delta P_L
\end{equation}
into $\hat{p}_L = \mathbb{E}[P_L]$, where $\mathbb{E}[\cdot]$ refers to the expected value, and $\Delta P_L$, which denotes the deviations from the expected value. Consequently, $\mathbb{E}[\Delta P_L]=0$ holds by definition. The same formulation can be found in \cite{appino2018use} or \cite{su2021optimal}. The \gls{pdf} $f_{P_L}$ is not assumed to follow any specific distribution (e.g., Gaussian), but needs to be available in a parametric form so that its integral over a specific domain can be computed. With that, the \gls{pdf} referring to the deviations from the expected value can be expressed as $f_{\Delta P_L}(z) = f_{P_L}(z+\hat{p}_L)$.

\subsection{Uncertainty Allocation}
\label{sec:uncertainty-splitting}

The main idea of the present work is the partitioning of the uncertain prosumption deviation $\Delta P_L$ into two parts, i.e., one part that is fed into the \gls{bess}, and one part that is injected into the grid. To this end, we introduce an upper and a lower bound of \gls{bess} charging power, i.e., $\ubar{x} \in \mathbb{R}^{-}_0$ and $\Bar{x} \in \mathbb{R}^{+}_0$. The prosumption uncertainty in the interval $[\ubar{x},\Bar{x}]$ is assigned to the \gls{bess}, and the remaining uncertainty is injected into the grid. We introduce two random variables that depend on $\Delta P_L$, i.e.,
\begin{subequations}    
\begin{align}
\label{eq:transformation-rv-b}
    \Delta P_{L\rightarrow B} &=
     \begin{cases}
       \ubar{x}  & \Delta P_{L} \le \ubar{x} \\
       \Delta P_{L} & \ubar{x} < \Delta P_{L} < \Bar{x}  \\
       \Bar{x} & \Bar{x} \le \Delta P_{L} \\ 
     \end{cases} \\
\label{eq:transformation-rv-g}
    & \nonumber \\
    \Delta P_{L\rightarrow G} &= 
     \begin{cases}
       \Delta P_{L} - \ubar{x}  & \Delta P_{L} \le \ubar{x} \\
       0 & \ubar{x} < \Delta P_{L} < \Bar{x}  \\
       \Delta P_{L} - \Bar{x} & \Bar{x} \le \Delta P_{L}. \\ 
     \end{cases} 
\end{align}
\end{subequations}
$\Delta P_{L\rightarrow B}$ and $\Delta P_{L\rightarrow G}$ represent the prosumption uncertainties that are shifted to the \gls{bess} or the grid, respectively.
With that, $\ubar{x}$ and $\Bar{x}$ reflect the minimal and maximal deviation from the expected prosumption $\hat{p}_L$ that is still fully compensated by the \gls{bess}. 
Note that $\ubar{x}$ and $\Bar{x}$ do not result from any forecasting process but can be chosen freely\footnote{To be precise, $\ubar{x}$ and $\Bar{x}$ can be chosen freely, but they are constrained, i.e., they are restricted by $\ubar{x}\le0$, $\Bar{x}\ge0$, and the respective \gls{bess}'s physical limits defined in \Cref{sec:battery-model}.} 
and thus be integrated into an optimization framework as decision variables, see \Cref{sec:optimization-problem}. 
For more information on the general transformation of random variables, see \cite{shankar2021probability} or \cite{edition2002probability}. 
Additionally, it should be noted that
\begin{equation}
\label{eq:uncertainty-distribution-definition}
    \Delta P_{L\rightarrow B} + \Delta P_{L\rightarrow G} = \Delta P_L
\end{equation}
holds by definition. This simply means that the sum of $\Delta P_{L\rightarrow B}$ and $\Delta P_{L\rightarrow G}$ accounts for the total prosumption uncertainty $\Delta P_L$.

\begin{figure*}[!htb]
    \centering
    \subfloat[Prosumption \gls{pdf}]{%
        \includegraphics[height=3cm]{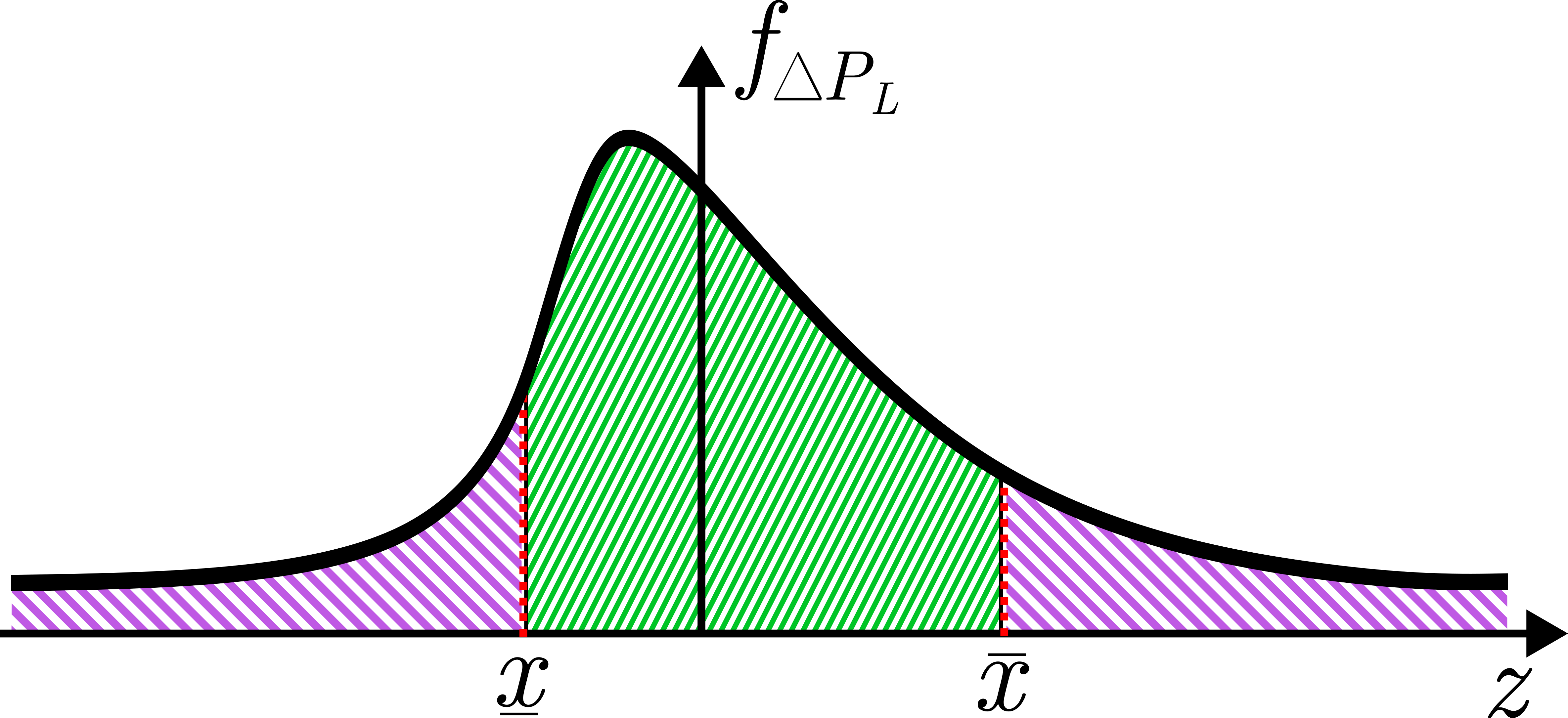}%
        \label{fig:pdf-prosumption}%
    }\hfill
    \subfloat[\gls{bess} Power \gls{pdf}]{%
        \includegraphics[height=3cm]{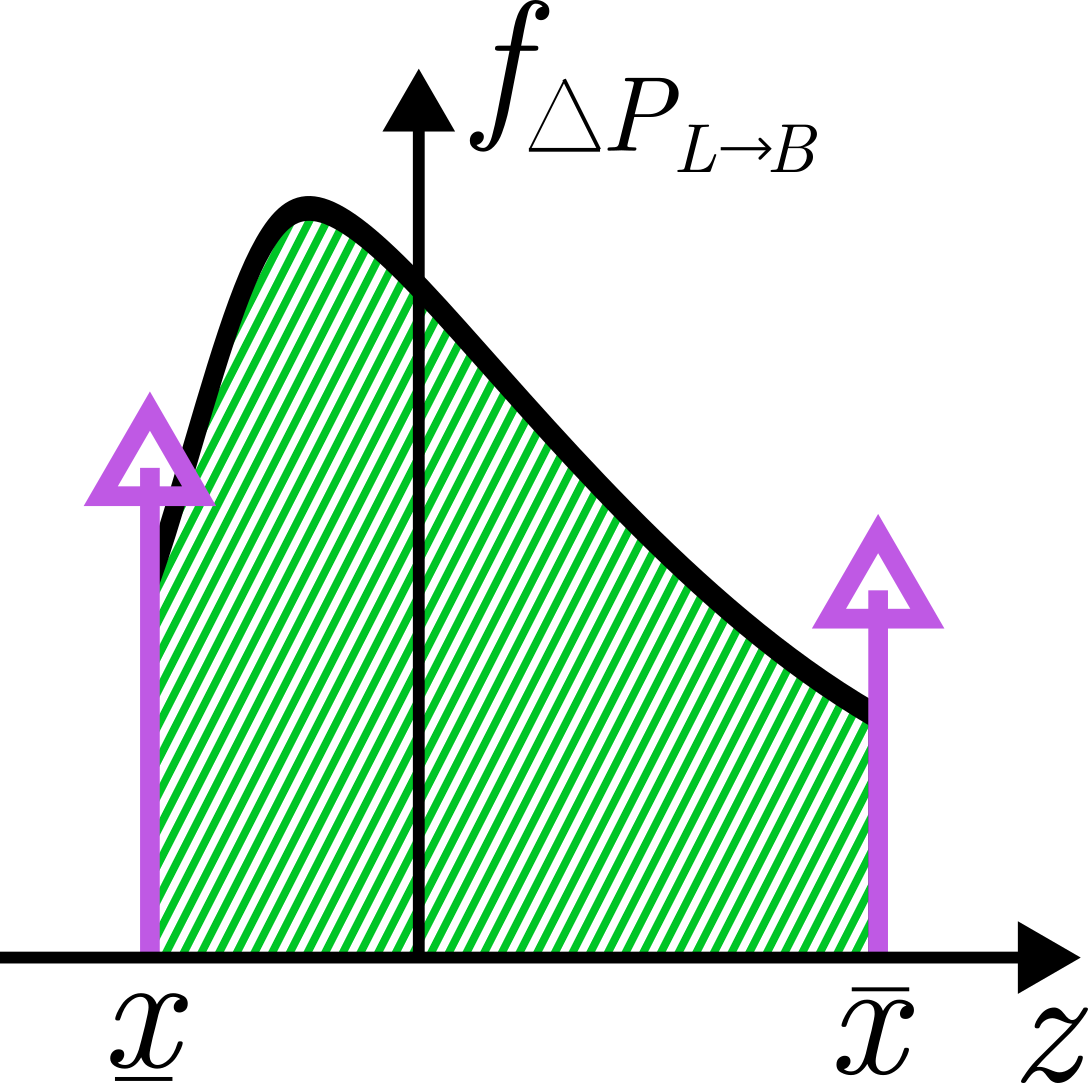}%
        \label{fig:pdf-battery}%
    }\hfill
    \subfloat[Grid Power \gls{pdf}] {%
        \includegraphics[height=3cm]{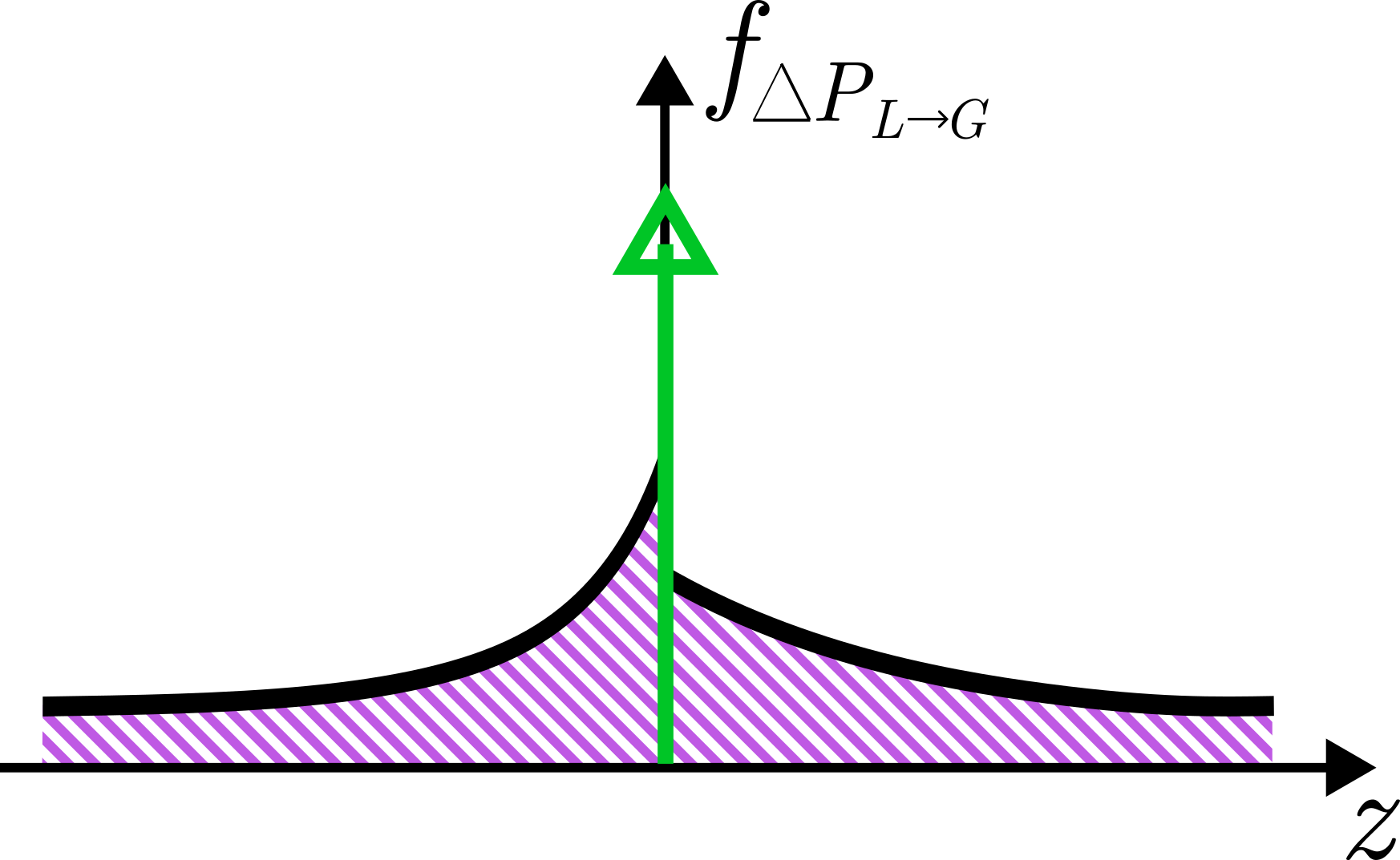}%
        \label{fig:pdf-grid}%
    }
    \caption{The continuous prosumption \gls{pdf} (a) can be represented by the two dependent distributions (b) and (c), containing continuous and discrete parts. The \gls{bess} power uncertainty (b) is defined over a closed interval. The grid power uncertainty (c) contains a discrete event at zero.}
    \label{fig:pdfs}
\end{figure*}

\autoref{fig:pdf-prosumption} illustrates the presented uncertainty allocation, with all uncertainty realizations within the interval $[\ubar{x}, \Bar{x}]$ being fully compensated by the \gls{bess} alone. This choice of compensation is advantageous since prosumption \glspl{pdf} typically have larger values around their expected values than in areas further away, essentially reducing the probability of compensation by the grid. For realizations outside of $[\ubar{x}, \Bar{x}]$, the occurring power discrepancy is provided by the \gls{bess} and the grid combined.

Furthermore, it should be noted that $\Delta P_{L\rightarrow B}$ and $\Delta P_{L\rightarrow G}$ contain both continuous and discrete parts, classifying them as so-called mixed random variables \cite{shankar2021probability}. The handling and interpretation of mixed random variables are not always apparent \cite{weld2017modeling}, but they have been proven useful to model different phenomena with mostly zero-inflated data such as air contamination \cite{owen1980estimation} or rainfall events \cite{li2013bivariate}. However, what sets our work apart from most research on mixed random variables is that the mixed distribution is intentionally designed, whereas in related studies, mixed distributions typically arise as natural phenomena from data. Rather than addressing the challenges posed by mixed distributions, this paper emphasizes the advantageous properties of mixed random variables and explores how they can be leveraged for modeling \gls{bess} scheduling.

\subsection{Computation of Mixed Distribution Functions and Expected Values}
\label{sec:mixedRandomVariables}
For a continuous random variable $X$ with \gls{pdf} $f_{X}$, the resulting distribution $f_Y$ of a mixed transformation $Y=g(X)$ can be derived by the superposition of continuous parts $f_{Y_C}$ and discrete parts $f_{Y_D}$ according to $f_Y(y) = f_{Y_C}(y) + f_{Y_D}(y)$. The two terms are derived and explained in \cite{shankar2021probability} and can be computed according to
\begin{subequations}
\begin{align}
        f_{Y_C}(y) &= \sum_{i=1}^{N} \left. \frac{f_{X}(x)}{\left| \frac{dy}{dx} \right|} \right|_{x_i=g^{-1}_{i}(y)} \label{eq:mixed-trafo-conti}\\
        &\nonumber\\
        f_{Y_D}(y) &= \sum_{i=1}^{M} \mathbb{P}(Y=y_{D,i}) \delta (y-y_{D,i}).
\end{align}
\end{subequations}
$N$ denotes the number of continuous parts of the mixed transformation, whereas $M$ refers to the number of discrete parts. In \Cref{eq:mixed-trafo-conti}, it is assumed that the inverses of the continuous parts of the transformation $g^{-1}_i$ exist.\footnote{In our case, the inverse functions are $g^{-1}_1(P_{L\rightarrow B}) = P_{L\rightarrow B}$ for the \gls{bess}, and $g^{-1}_1(P_{L\rightarrow G}) = P_{L\rightarrow G} + \ubar{x}$ and $g^{-1}_2(P_{L\rightarrow G}) = P_{L\rightarrow G} + \Bar{x}$ for the grid.}
$y_{D,i}$ are discrete events, and $\delta$ is the Dirac delta function.
With that the \glspl{pdf}\footnote{When dealing with mixed random variables, the term \textit{probability density function} is sometimes used to describe the entire distribution, even though this is not mathematically precise. A \gls{pdf} applies to continuous components only, while a \gls{pmf} refers to discrete ones. Since no standard term covers both aspects of mixed distributions, the term \gls{pdf} is used in the present paper for simplicity.} of $\Delta P_{L\rightarrow B}$ and $\Delta P_{L\rightarrow G}$ can be computed using the unit step function $U(z)$:
\begin{subequations}
\begin{align}
\begin{split}
    f_{\Delta P_{L\rightarrow B}}(z) =\ &p_1 \delta(z-\ubar{x}) + p_2 \delta(z-\Bar{x}) \\
    &+ f_{\Delta P_L}(z) [U(z-\ubar{x}) - U(z-\Bar{x})]
\end{split}
    & \\
    & \nonumber \\
\begin{split}
    f_{\Delta P_{L\rightarrow G}}(z) =\ &(1 - p_1 - p_2) \delta(z) + f_{\Delta P_L}(z + \Bar{x}) U(z) \\
    &+ f_{\Delta P_L}(z + \ubar{x}) [U(z + \infty) - U(z)],
\end{split}
\end{align}
\end{subequations}
where the probabilities $p_1 = \mathbb{P}(\Delta P_{L\rightarrow B} = \ubar{x})$ and $p_2 =\mathbb{P}(\Delta P_{L\rightarrow B} = \Bar{x})$ correspond to the discrete events. The unit step function $U(z)$ restricts the domains in which the continuous parts are defined.
The resulting \gls{pdf} of the \gls{bess} power uncertainty is illustrated in \autoref{fig:pdf-battery}. The uncertainties of the original prosumption outside the interval $[\ubar{x}, \Bar{x}]$ are mapped to the discrete events $\ubar{x}$ or $\Bar{x}$ respectively, effectively limiting the maximal and minimal \gls{bess} power deviation to a fixed interval. The other portion of the prosumption uncertainty remains unchanged. \autoref{fig:pdf-grid} displays the uncertainties shifted to the grid as a zero-inflated distribution; all uncertainties within $[\ubar{x}, \Bar{x}]$ are mapped to zero, whereas the uncertainty tails are shifted horizontally, but remain otherwise untouched. It can easily be seen that all displayed distributions are valid as they are all non-negative and integrate to one (counting the discrete and continuous events).


Furthermore, due to the independent choice of $\ubar{x}$ and $\Bar{x}$, the uncertainties can be asymmetrically allocated between the grid and the \gls{bess}.
This property is shown to be useful in \Cref{sec:model-application}, but also implies that the deviations $\Delta P_{L\rightarrow B}$ and $\Delta P_{L\rightarrow G}$ can have non-zero expected values, i.e.,
\begin{subequations}
\begin{align}
    \label{eq:exp-battery}
    \mathbb{E}[\Delta P_{L\rightarrow B}] &= p_1\ubar{x} + \int_{\ubar{x}}^{\Bar{x}} z f_{\Delta P_L}(z)\, dz  + p_2\Bar{x}\\
    & \nonumber \\
    \label{eq:exp-grid}
    \begin{split}
    \mathbb{E}[\Delta P_{L\rightarrow G}] &= \int_{-\infty}^{0} z f_{\Delta P_L}(z + \ubar{x})\, dz  \\
    &\ + \int_{0}^{\infty} z f_{\Delta P_L}(z + \Bar{x})\, dz.
    \end{split}
\end{align}
\end{subequations}
However, because the expected value of a sum of random variables is equal to the sum of their individual expected values, regardless of whether they are independent \cite{shankar2021probability}, 
\begin{align}
\begin{split}
        \mathbb{E}[\Delta P_L] &= \mathbb{E}[\Delta P_{L\rightarrow B} + \Delta P_{L\rightarrow G}]  \\
    &= \mathbb{E}[\Delta P_{L\rightarrow B}] + \mathbb{E}[\Delta P_{L\rightarrow G}] = 0
\end{split}
\end{align}
is valid and can be either used as an alternative to \Cref{eq:exp-battery} or \Cref{eq:exp-grid}, or simply for validation purposes. 

\subsection{Battery State Evolution}
\label{sec:battery-model}

In this section, the \gls{bess} model with its state evolution is derived. We describe the state evolution in discrete time with time variable $k\in \mathcal{K}\subset\mathbb{N}_0$.
Because of the described probabilistic approach with the allocation of uncertainties in \Cref{sec:uncertainty-splitting}, the \gls{bess} energy state becomes a random variable $E$, which we split into two parts, i.e.,
\begin{equation}
    E(k) = e(k) + \Delta E(k).
\end{equation}
$e(k) \in \mathbb{R}$ is introduced as the nominal battery state and reflects the energy state based on the previous expected prosumptions $\hat p_L$, whereas the probabilistic battery state $\Delta E(k)$ is a random variable that depends on the previous probabilistic prosumptions $\Delta P_L$.
To be concise,
\begin{align}
    e(k) &= h_1(\hat{p}_L(k), \hat{p}_L(k-1), ..., \hat{p}_L(0)) \nonumber \\
    \Delta E(k) &= h_2(\Delta P_L(k), \Delta P_L(k-1), ..., \Delta P_L(0)). \nonumber
\end{align}
It is important to note, that even though the nominal battery state $e$ corresponds to the expected prosumption, it does not reflect the expected battery state, hence the different name. This phenomenon arises due to the handling of asymmetric prosumption \glspl{pdf} as well as the possibility of asymmetric uncertainty shifts. The nominal battery state $e$ coincides with the expected battery state $\mathbb{E}[E]$, only if \Cref{eq:exp-battery} permanently equals zero. This behavior is exemplarily analyzed in \Cref{sec:results}. 

The \gls{bess} power can be described with a similar approach via
\begin{equation}
\label{eq:bat-power-split}
    P_B(k) = p_B(k) + \Delta P_{L\rightarrow B}(k)
\end{equation}
with $p_B$ being the nominal battery power. 

The evolution of the \gls{bess} state can be formulated:
\begin{align}
    E(k+1) &= E(k) - t P_B(k) - t \mu \left | P_B(k) \right | \label{eq:general-battery-evo} \\
    &= e(k) + \Delta E(k) - t p_B(k) - t \Delta P_{L\rightarrow B}(k) \nonumber \\ &\ \ \ \ \ \ \ \ \  - t \mu \left | p_B(k) + \Delta P_{L\rightarrow B}(k) \right |. \nonumber
\end{align}
The variable $\mu$ is a loss coefficient that models charging and discharging losses, $t$ denotes the time between the discrete time points $k$ and $k+1$. 
To separate evolution of the nominal and probabilistic battery states, $E$ is approximated by using the triangle inequality $|a+b| \le |a| + |b|$:
\begin{align}
    E(k+1)&= e(k) + \Delta E(k) - t p_B(k) - t \Delta P_{L\rightarrow B}(k) \nonumber \\ &\ \ \ \ \ \ \ \ \  - t \mu \left | p_B(k) + \Delta P_{L\rightarrow B}(k) \right | \nonumber \\
    \begin{split}
    &\ge e(k) + \Delta E(k) - t p_B(k) - t \Delta P_{L\rightarrow B}(k) \\ &\ \ \ \ \ \ \ \ \  - t \mu \left | p_B(k) \right | - t \mu \left |\Delta P_{L\rightarrow B}(k) \right |.
\end{split}
\end{align}
This is justifiable because only comparatively small portions of the total energy state are approximated, i.e., the \gls{bess}'s losses. Additionally, this approximation can be considered conservative since the total energy state is underestimated. 

Hence, the nominal battery state evolution
\begin{equation}
    e(k+1) = e(k) - t p_B(k) - t\mu \left | p_B(k) \right |
\end{equation}
and the probabilistic evolution 
\begin{equation}
    \label{eq:bat-evo-probabilistic}
    \Delta E(k+1) = \Delta E(k) - t \Delta P_{L\rightarrow B}(k) - t \mu  \left | \Delta P_{L\rightarrow B}(k) \right |
\end{equation}
can be formulated separately. 
\Cref{eq:bat-evo-probabilistic} consists of the sum of three random variables, making its computation challenging. If one neglects stochastic losses $\mu  \left | \Delta P_{L\rightarrow B}(k) \right | \approx 0$ and assumes independence of prosumption between subsequent time increments $\Delta P_L(k)$ and $\Delta P_L(k+1) \ \forall k \in \mathcal{K}$, the probabilistic battery state can be computed by convolution. However, such an assumption is difficult to justify: For instance, even a simple unforeseen cloud formation can affect a \gls{pv}'s power production over several time increments, an effect that would be entirely overlooked under this assumption. 
Consequently, some \gls{bess} operation models are specifically designed to account for this dependence through sequential decision-making processes \cite{gross2020stochastic, pacaud2024optimization}.
An alternative method to compute \Cref{eq:bat-evo-probabilistic} would involve finding a suitable approximation of the joint probability function. In \cite{austnes2023probabilistic}, a way to approximate a joint dependence for probabilistic load forecasting using empirical copulas is presented. 
However, applying such an approach to the problem at hand is challenging due to the complex dependence between probabilistic battery power and probabilistic battery state, where the latter essentially encapsulates all prior uncertainties of probabilistic battery power. This complexity is further compounded by the fact that the probabilistic battery state is influenced by decision variables $\ubar{x}(k)$ and $\Bar{x}(k)$, which are supposed to be found by solving a subsequent optimization problem.

Although a suitable solution for solving \Cref{eq:bat-evo-probabilistic} has not yet been identified, the equation remains valuable for our approach: The restriction of \gls{bess} power uncertainties to the fixed interval $[\ubar{x}(k), \Bar{x}(k)]$ (see \autoref{fig:pdf-battery}) implies that the resulting energy uncertainty accumulated in the \gls{bess} is also bounded. This insight can be used to calculate the evolution of the minimum and maximum battery state based on \Cref{eq:bat-evo-probabilistic}:
\begin{subequations}
\begin{align}
    \Delta E_{\min}(k+1) &= \Delta E_{\min}(k) - t\Bar{x}(k) - t\mu\Bar{x}(k) \\
    \Delta E_{\max}(k+1) &= \Delta E_{\max}(k) - t\ubar{x}(k) + t\mu\ubar{x}(k),
\end{align}
\end{subequations}
with $\ubar{x}(k)\le 0$ and $\Bar{x}(k) \ge 0$. 
Therefore, we can compute the minimum and maximum battery states, the nominal battery state, and the expected battery state (which can be derived from \Cref{eq:exp-battery}).
Nevertheless, the lack of a full distribution function of the \gls{bess} state is a limitation of the model depicted. 
This restricts the ability to shift \gls{bess} uncertainties to the grid, meaning the uncertainty bandwidth within the \gls{bess} is non-decreasing. 

Typically, the \gls{bess} energy and power constraints are expressed as 
\begin{subequations}
\begin{align}
    E(k) &\stackrel{!}{\in} [\ubar{e}, \Bar{e}] \\
    P_B(k) & \stackrel{!}{\in} [\ubar{p}_B, \Bar{p}_B],
\end{align}
\end{subequations}
with $\ubar{e}$, $\Bar{e}$ and $\ubar{p}_B$, $\Bar{p}_B$ being the minimum and maximum \gls{bess} capacity and power, respectively.
In our approach, the \gls{bess} constraints can be reformulated to
\begin{subequations}
\begin{align}
    e(k) + \Delta E_{\min}(k) &\ge \ubar{e} \\
    e(k) + \Delta E_{\max}(k) &\le \Bar{e} \\
    p_B(k) + \ubar{x}(k) &\ge \ubar{p}_B \\
    p_B(k) + \Bar{x}(k) &\le \Bar{p}_B .
\end{align}
\end{subequations}
Their compliance can be ensured at all times, implying that the physical constraints of the \gls{bess} are respected.

\subsection{Grid Exchange}
\label{sec:grid-exchange}

The grid power
\begin{equation}
\label{eq:power-grid-split}
    P_G = p_G + \Delta P_{L\rightarrow G}
\end{equation}
is divided into two components - nominal and probabilistic - mirroring the partitioning of the \gls{bess} power in \Cref{eq:bat-power-split}. Based on the lossless power connection in \Cref{eq:power-balance}, the prosumption uncertainty allocation in \Cref{eq:uncertainty-distribution-definition} and the partitioning of random variables in two parts in \Cref{eq:prosumption}, (\ref{eq:bat-power-split}) and (\ref{eq:power-grid-split}), it directly follows that
\begin{equation}
    \hat{p}_L = p_B + p_G.
\end{equation}
It can be seen that even though the nominal powers depend solely on the expected prosumption $\hat{p}_L$, they do not reflect the expected powers. This discrepancy arises because the additional expected values of $\Delta P_{L\rightarrow B}$ and $\Delta P_{L\rightarrow G}$ (see \Cref{eq:exp-battery} and (\ref{eq:exp-grid})) cancel each other out during the derivation.

In the following, $p_G$ represents the \acrfull{ds}, while $\Delta P_{L\rightarrow G}$ describes the deviations from the \gls{ds}. These deviations should be minimized, ideally being zero. 
With the chosen structure of $\Delta P_{L\rightarrow G}$ to specifically include zero-inflated events, deviations of zero can appear with non-zero probability (see \autoref{fig:pdf-grid}).

\section{Optimization Problem}
\label{sec:optimization-problem}

In this section, the derived model is summarized within the framework of an optimization problem. The goal is to determine an optimized probabilistic day-ahead \gls{bess} schedule and to establish an optimized probabilistic \gls{ds} with the grid. The problem can be formulated as follows:
\begin{subequations}
\begin{align}
    \underset{p_B, \ubar{x}, \Bar{x}}{\min} \ C(p_G^+,\ p_G^-, \ & p_B^+,\  p_B^-, \ p_1,\ p_2,\ \mathbb{E}[\Delta P_{G}^{<0}],\ \mathbb{E}[\Delta P_{G}^{>0}]) \nonumber \\
    \text{s.t. for all } k \in \mathcal{K} \nonumber \\
     e(k+1) = e(&k) - t p_B(k) - t \mu p_B^+(k) + t \mu p_B^-(k) \label{op:e_nom} \\
    \Delta E_{\min}(k+1) &= \Delta E_{\min}(k) - t\Bar{x}(k) - t \mu \Bar{x}(k)  \label{op:e_min}\\
    \Delta E_{\max}(k+1) &= \Delta E_{\max}(k) - t\ubar{x}(k) + t \mu \ubar{x}(k)\label{op:e_max} \\
    %
    & \nonumber \\
    e(k) &+ \Delta E_{\min}(k) \ge \ubar{e} \label{op:e_constr1}\\
    e(k) &+ \Delta E_{\max}(k) \le \Bar{e} \label{op:e_constr2}\\
    p_B(k) &+ t\ubar{x}(k) \ge \ubar{p}_B \label{op:e_constr3} \\
    p_B(k) &+ t\Bar{x}(k) \le \Bar{p}_B \label{op:e_constr4}\\
    & \nonumber \\
    p_G(k) &= \hat{p}_L(k) - p_B(k)\label{op:pg_nom} \\
    \mathbb{E}[\Delta P_{G}^{<0}(k)] &= \int_{-\infty}^{0} z f_{k, \Delta P_L}(z + \ubar{x}(k))\, dz \label{op:pg_exp_low}\\
    \mathbb{E}[\Delta P_{G}^{>0}(k)] &= \int_{0}^{\infty} z f_{k, \Delta P_L}(z + \Bar{x}(k))\, dz  \label{op:pg_exp_high}\\
    p_1(k) &= \int_{-\infty}^{\ubar{x}(k)} f_{k, \Delta P_L}(z)\, dz \label{op:p1} \\
    p_2(k) &= 1 - \int_{-\infty}^{\Bar{x}(k)} f_{k, \Delta P_L}(z)\, dz \label{op:p2} \\
    & \nonumber \\
    p_{G}(k) &= p_G^+(k) + p_G^-(k) \label{op:pg_split1} \\
     -p_G^+(k) \cdot p_G^-(k) &\le 10^{-8} \label{op:pg_split2} \\
     p_G^+(k) \ge\ &0,\  \ p_G^-(k) \le 0 \label{op:pg_split3} \\
    & \nonumber \\
    p_{B}(k) &= p_B^+(k) + p_B^-(k)  \\
     -p_B^+(k) \cdot p_B^-(k) &\le 10^{-8}  \\
     p_B^+(k) \ge\ &0,\  \ p_B^-(k) \le 0  \\
     & \nonumber \\
    \ubar{x}(k) &\le 0  \\
    \Bar{x}(k) &\ge 0. 
    %
\end{align}
\end{subequations}
The decision variables of the problem are $p_B(k)$, $\ubar{x}(k)$ and $\Bar{x}(k)$. The input variables, i.e., $\hat{p}_L(k)$ and $f_{k, \Delta P_L}(z)$, are the result of a preceding probabilistic forecasting process.

Equations (\ref{op:e_nom}-\ref{op:e_max}) refer to the \gls{bess} state evolutions. Instead of having fixed scheduling set points, the charging power and the energy state are characterized by respective bounds. Equations (\ref{op:e_constr1}-\ref{op:e_constr4}) ensure compliance with the \gls{bess}'s physical limits.
\Cref{op:pg_nom} computes the nominal grid power, which is used as the basis for the \gls{ds}. Thus, deviations from $p_G$ must be quantified, which is described in Equations (\ref{op:pg_exp_low}-\ref{op:p2}). To elaborate, incorporating $p_1$ and $p_2$ into the cost function minimizes the likelihood of deviations from the computed \gls{ds}. Specifically, $p_2$ accounts for the probability of requiring more power, while $p_1$ addresses the probability of needing less. Additionally, the extent of these deviations can be captured through $\mathbb{E}[\Delta P_{G}^{<0}]$ and $\mathbb{E}[\Delta P_{G}^{>0}]$, which represent the expected magnitude of downward or upward deviations, respectively, as defined by \Cref{eq:exp-grid}.
The respective integrals depend on the decision variables and thus need to be computed dynamically during the optimization process. For that, the integrals are approximated using numerical integration, specifically, Simpson's 1/3 rule \cite{brass2011quadrature}.

The nominal grid power $p_G$ is split into two mutually exclusive components in Equations (\ref{op:pg_split1}-\ref{op:pg_split3}): a positive $p_G^+$ and a negative $p_G^-$. This enables the cost function to define electricity purchase and sale prices separately, making it possible to model more realistic market scenarios. The nominal battery power $p_B$ is divided similarly to accurately account for energy losses.

To summarize, the presented method is an expectation-based nonlinear stochastic optimization problem that is suitable for typical \gls{bess} applications, such as maximizing self-sufficiency, as well as actively reducing quantified grid uncertainties.
Thereby, the key innovation lies in the analytical and coherent treatment of uncertainties described by mixed random variables. Unlike previous approaches, where uncertainties are either entirely absorbed by the grid \cite{su2021optimal} or fully shifted into the \gls{bess} \cite{appino2018use}, the presented method enables a quantified and somewhat balanced distribution of uncertainties between the \gls{bess} and the grid. This allows the optimization problem to simultaneously compute an optimized probabilistic grid and probabilistic \gls{bess} schedule --- a capability that, to the best of the authors' knowledge, has not yet been published and analyzed in the literature.
This versatility opens up a wide range of potential applications, several of which are analyzed in the following section.

\section{Model Application}
\label{sec:model-application}

This section presents the application of the proposed model, highlighting different features in three distinct cases. Based on real-world data, we compute probabilistic forecasts for the prosumption and solve the day-ahead optimization problem for three different cost functions. In Case 1, the model aims to reduce electricity costs by increasing the household’s self-sufficiency. Case 2 extends Case 1 by additionally minimizing prosumption uncertainties within the grid. Case 3 combines cost reduction with minimizing grid uncertainties during a critical timeframe.

The model is implemented in Python using the optimization modeling language Pyomo \cite{bynum2021pyomo} with IPOPT \cite{wachter2006implementation} as the solver.\footnote{The implementation is published as an open-access repository along with this paper \cite{Pinter_Day-Ahead_Battery_Scheduling_2024}. All presented results can be reproduced using the provided code.} All optimization processes are solved on a laptop with an Intel Core i7-1355U CPU at 1.70 GHz and 32 GB RAM in less than one minute.

\subsection{Experimental Setup}
\label{sec:experimental-setup}

We use the openly available power data from the Open Power System Data platform \cite{opendataset} of a single residential house referenced in the dataset as ‘residential4’ located in southern Germany. This specific building is characterized by a high \gls{pv} output that exceeds its consumption most of the time. Therefore, the net prosumption is mostly negative. The available data covers the years 2016 to 2018 and serves as the basis for creating the probabilistic forecasts. For forecasting, the Temporal Fusion Transformer \cite{lim2021temporal} using the DARTS framework \cite{JMLR:v23:21-1177} is used. The forecasting model predicts 24 hours of prosumption in an hourly resolution, starting at 06:00 each day. Witherating 100 forecasts per hour, respective quantiles from  1\% to 99\% in equidistant increments are estimated. These forecasts are based on the previous seven days of data, along with additional covariates such as year, month, weekday, and hour of the day, covering both the historical data and the forecasted period. 
Note that the forecast model does not include weather data, as they are not part of the available dataset. While the inclusion of data, such as global radiation, would likely reduce the forecast uncertainty concerning the \gls{pv} production, the forecast uncertainty in consumption would remain the same.
The first year of data is used for training the forecasting model and the following 90 days are used for validation and early stopping. The forecasted quantiles are smoothed.

\begin{figure}
    \centering
    \includegraphics[width=1.0\linewidth]{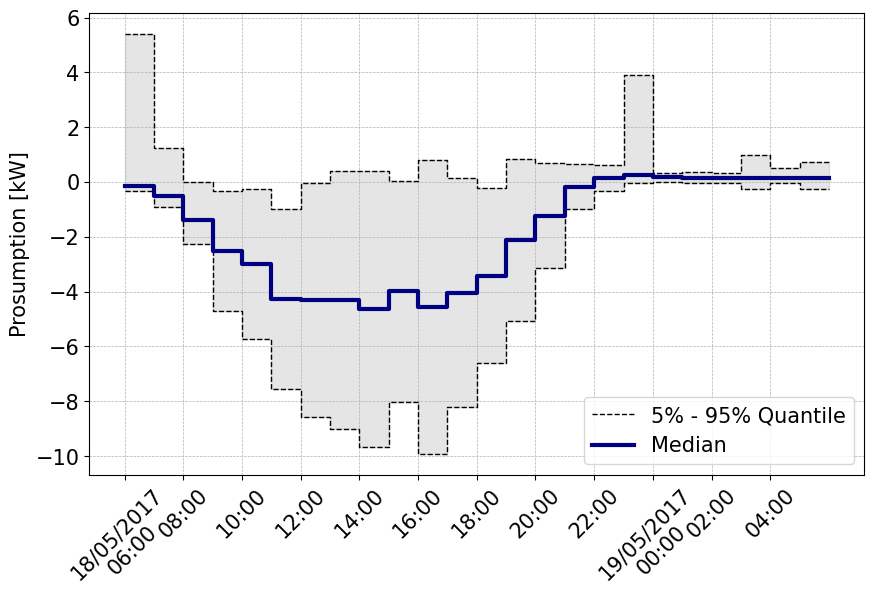}
    \caption{Forecasted residential prosumption $P_L$ for a typical day. The \gls{pv} production is expected to dominate the consumption most of the time.}
    \label{fig:forecast-prosumption}
\end{figure}

\begin{figure}[h]
    \centering
    \includegraphics[width=1.0\linewidth]{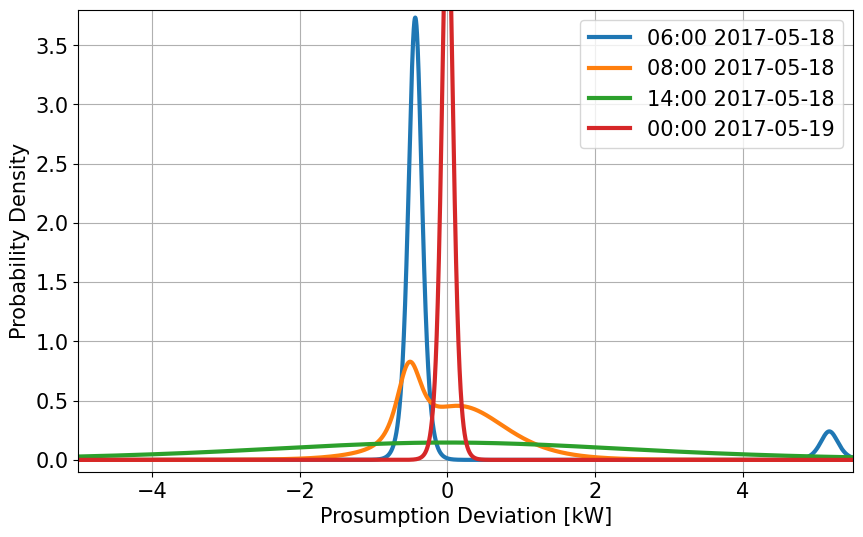}
    \caption{The overall shape, variance, and level of asymmetry of the forecasted prosumption \glspl{pdf} $f_{\Delta P_L}$ vary highly throughout the day. All \glspl{pdf} have an expected value of zero.}
    \label{fig:pdf-prosumption-forecast}
\end{figure}

In order to incorporate the forecasts into the optimization framework, they must be available in parametric form. In the present study, the uncertainties are represented by random variables whose \glspl{cdf} are assumed to follow the sum of two logistic functions:
\begin{equation}
\label{eq:cdf-formula}
    F_{P_L}(z) =  \frac{\omega_1}{1+\exp({-\omega_2 (z - \omega_3)})} + \frac{\omega_4}{1+\exp({-\omega_5 (z - \omega_6)})}.
\end{equation}
The weights $w_i, i \in [1,6]$ are obtained via a fitting process to the forecasted quantiles. This form is selected for its ability to capture the asymmetry often observed in real prosumption data. \autoref{fig:forecast-prosumption} illustrates the forecasted prosumption for a typical day of the selected building.
The respective \glspl{pdf} vary over time in their level of asymmetry, shape, and total variance, see \autoref{fig:pdf-prosumption-forecast}. 
At 06:00, a peak at 5 kW in the respective \gls{pdf} can be seen. This peak represents a high consumption with a low probability, potentially stemming from occasional residential early activities (e.g., using the washing machine, dishwasher, stove, ...).

The cost function of the optimization problem is chosen so that the total power exchange with the grid $p_G^2$ and the impact and probability of deviations from the \gls{ds} can be minimized, i.e.,
\begin{align}
\begin{split}
    C = \sum_{k = 0}^{K} c_1 {p_{G, k}^+}^2 + c_2 {p_{G, k}^-}^2 &+ c_3 p_{1, k} \mathbb{E}[\Delta P_{G, k}^{>0}] \\
    + c_4 p_{2, k} \mathbb{E}[\Delta P_{G, k}^{<0}].
\end{split}
\end{align}
The probabilistic forecasts are fed into three optimization problems, distinguished by the cost-function weights specified in \autoref{tab:objective-weights}. 
In Case 1, the only objective is to minimize the homeowner’s electricity costs by enhancing self-sufficiency. The weights $c_3$ and $c_4$ are set to zero, indicating that deviations from the \gls{ds} are not penalized. As a result, the optimization neither penalizes grid uncertainties nor uncertainties inside the \gls{bess}. Consequently, $\ubar{x}$ and $\Bar{x}$ are dispensable decision variables and are therefore set to zero. This formulation leads to a deterministic \gls{bess} schedule, with all uncertainties shifted toward the grid.
For Case 2, the primary objective remains to minimize electricity costs. Additionally, deviations from the \gls{ds} are slightly penalized.
Lastly, Case 3 seeks to balance minimizing electricity costs with reducing downward grid deviations during a critical timeframe. This case reflects a foreseeable future event, where it is unfavorable to have high downward deviations.
For example, suppose a power line is operating near full capacity and is at risk of overheating due to expected high \gls{pv} generation peaks \cite{ccakmak2022using}, it may make sense to use a \gls{bess} to partially compensate for even higher potential \gls{pv} generation peaks, thereby relieving the grid during this critical period. For all three cases, standard \gls{bess} parameters \cite{teslapowerwall} are used, see \autoref{tab:params}.

\begin{table}[ht]
    \centering
    \caption{Selected cost-function weights.}
    \begin{tabular}{l c c c c} 
        \toprule
        \textbf{Test Case} & $c_1$ & $c_2$ & $c_3$ & $c_4$  \\
        \midrule
        \textbf{Case 1}  & 2 & 1 & 0 & 0 \\
        \textbf{Case 2}  & 2 & 1 & 0.5 & 0.5 \\
        \multirow{2}{*}{\textbf{Case 3}}  & \multirow{2}{*}{2} & \multirow{2}{*}{1} & \multirow{2}{*}{2} &100 for $k \in [6,7]$,\\ &&&&2 otherwise\\
        \bottomrule
    \end{tabular}
    \label{tab:objective-weights}
\end{table}

\begin{table}[]
    \centering
    \caption{Selected \gls{bess} specifications.}
    \label{tab:params}
    \begin{tabular}{c c c c c}
        \toprule
        $\ubar{e}$ [kWh] & $\Bar{e}$ [kWh] & $\ubar{p}_B$ [kW] &  $\Bar{p}_B$ [kW] & $\mu$ [\%]\\
        \midrule
        0.0 & 13.5 & -5.0 &  5.0 & 5.0\\
        \bottomrule
    \end{tabular}
\end{table}

\subsection{Results}
\label{sec:results}

\begin{figure*}[!htb]
    \centering
    \subfloat{%
        \begin{minipage}[b]{0.45\textwidth}
            \centering
            \includegraphics[width=\textwidth]{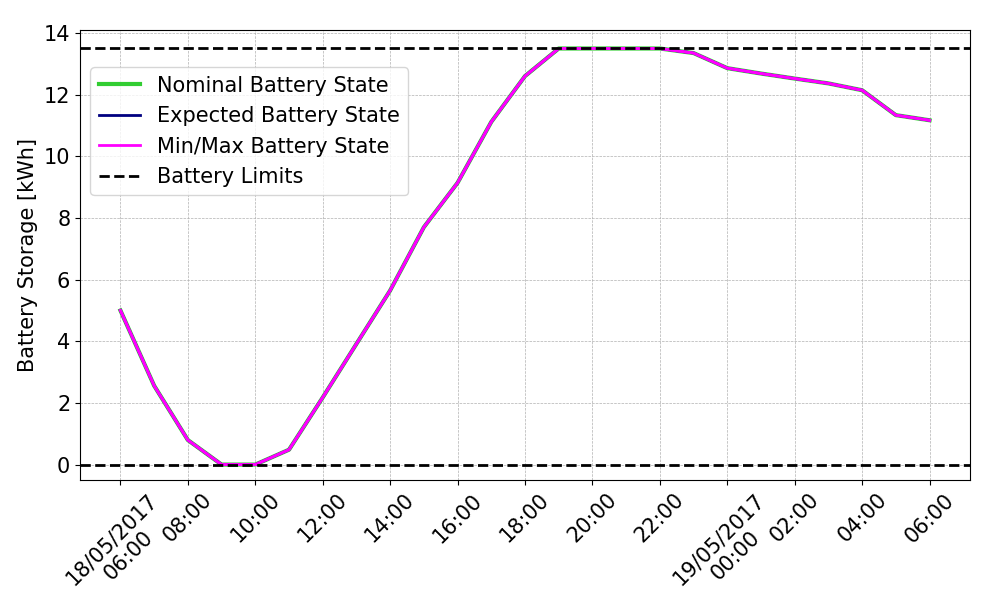}
            (a) Deterministic \gls{bess} Schedule
        \end{minipage}
    }\hfill
    \subfloat{%
        \begin{minipage}[b]{0.45\textwidth}
            \centering
            \includegraphics[width=\textwidth]{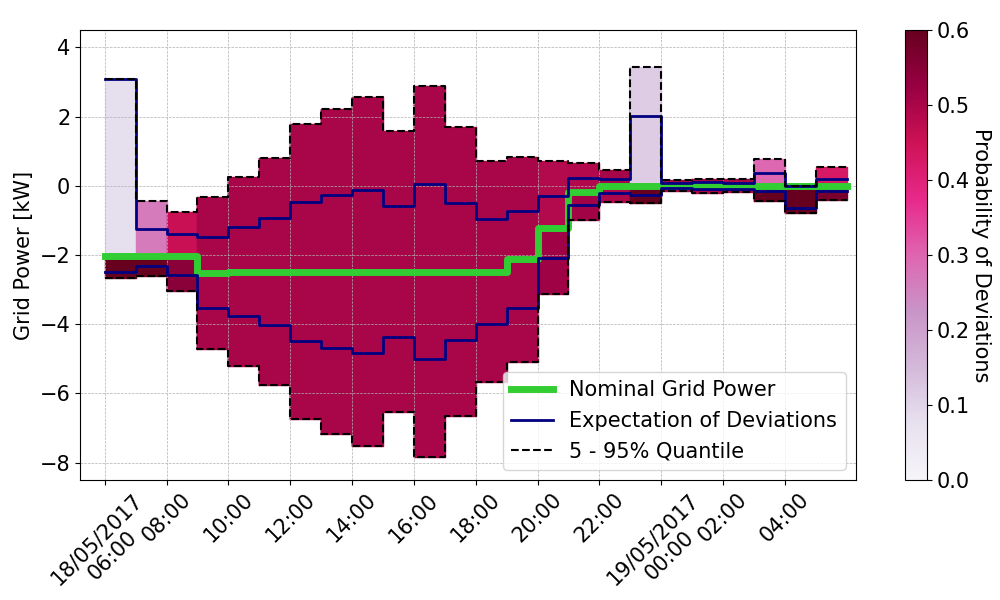}
            (b) Probabilistic Dispatch Schedule
        \end{minipage}
    }\\
    \caption{Case 1 - Minimizing residential electricity costs by enhancing self-sufficiency. All \gls{bess} states coincide, making the \gls{bess} schedule deterministic. All uncertainties are shifted toward the grid.}
    \label{fig:case1}
    
    \vspace{1em}
    
    \subfloat{%
        \begin{minipage}[b]{0.45\textwidth}
            \centering
            \includegraphics[width=\textwidth]{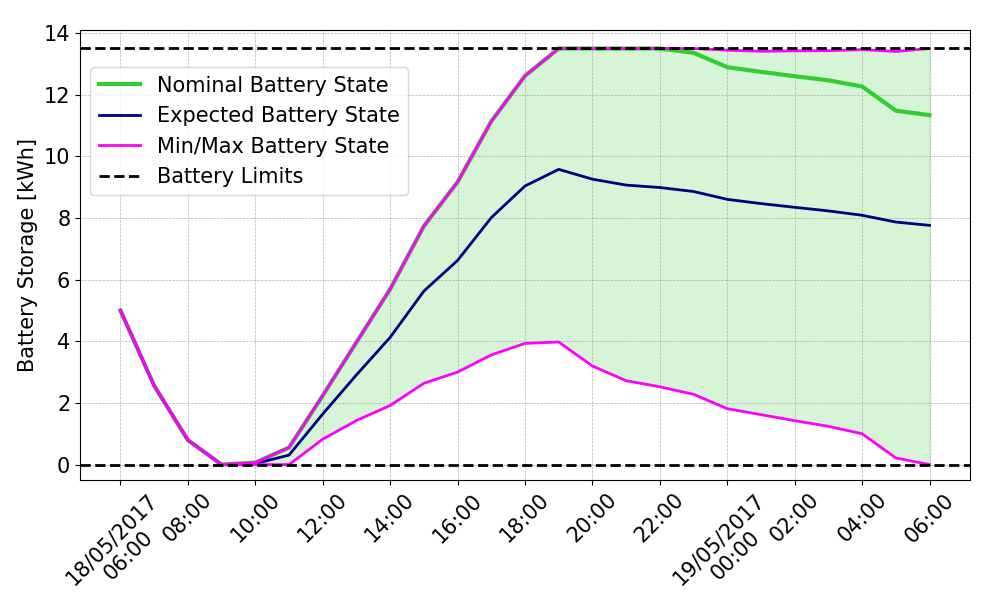}
            (c) Probabilistic \gls{bess} Schedule
        \end{minipage}
    }\hfill
    \subfloat{%
        \begin{minipage}[b]{0.45\textwidth}
            \centering
            \includegraphics[width=\textwidth]{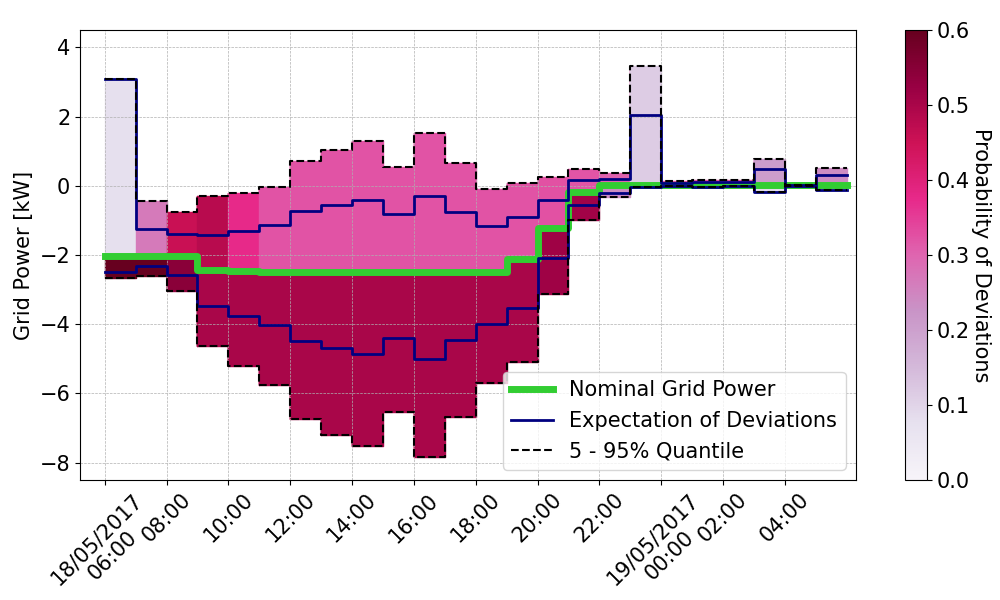}
            (d) Probabilistic Dispatch Schedule
        \end{minipage}
    }\\
    \caption{Case 2 - Minimizing residential electricity costs, while relieving the grid. The \gls{bess} reduces grid uncertainties without jeopardizing its main objective of minimizing residential electricity costs by enhancing self-sufficiency. Note that the nominal grid power of Case 1 and Case 2 align, indicating that both cases commit to a schedule with the same electricity costs.}
    \label{fig:case2}

    \vspace{1em}

    \subfloat{%
        \begin{minipage}[b]{0.45\textwidth}
            \centering
            \includegraphics[width=\textwidth]{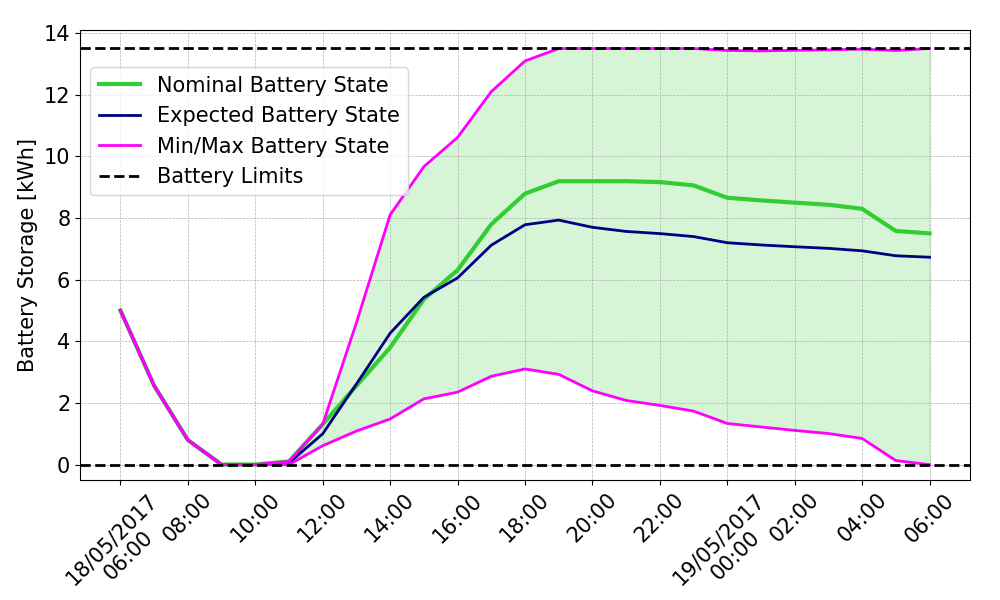}
            (e) Probabilistic \gls{bess} Schedule
        \end{minipage}
    }\hfill
    \subfloat{%
        \begin{minipage}[b]{0.45\textwidth}
            \centering
            \includegraphics[width=\textwidth]{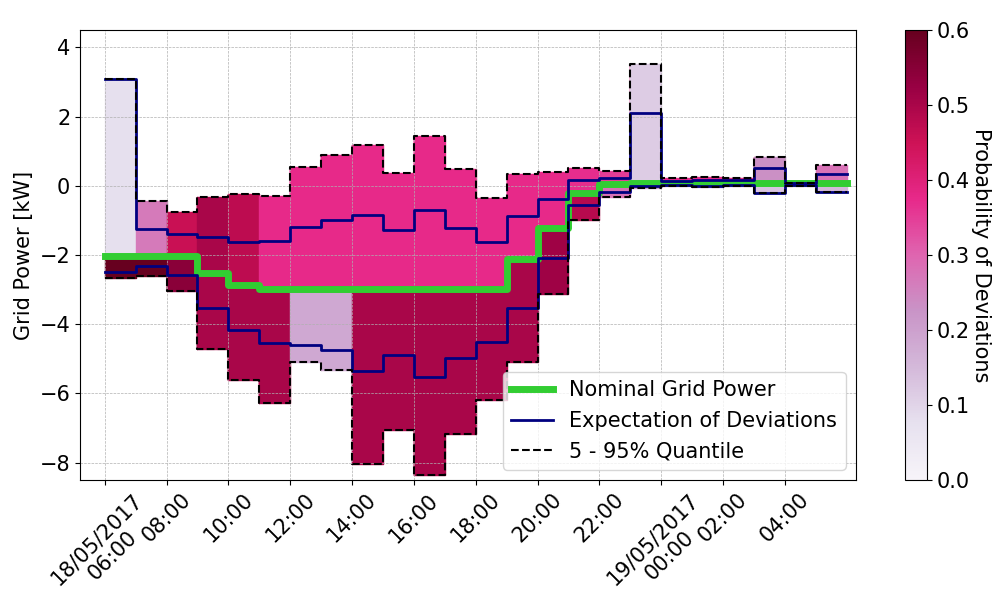}
            (f) Probabilistic Dispatch Schedule
        \end{minipage}
    }\\
    \caption{Case 3 - Reducing grid uncertainties during a critical timeframe. Between 12:00 and 14:00, the \gls{bess} is actively used to minimize the probabilities of downward deviations from the nominal grid power.}
    \label{fig:case3}
\end{figure*}

The results of the three introduced cases are presented in \autoref{fig:case1}, \autoref{fig:case2}, and \autoref{fig:case3}. The figures on the left-hand side display the \gls{bess} schedules, while the figures on the right-hand side illustrate the corresponding probabilistic \glspl{ds}. The \gls{bess} schedules are characterized by four different \gls{bess} states, as introduced in \Cref{sec:battery-model}.

The probabilistic \glspl{ds} provide the nominal grid trajectory, see \Cref{sec:grid-exchange}. In addition, they visualize the probabilities of deviations from the nominal grid power, with downward and upward deviations shown separately. To capture the plausible range of these deviations, the 5\% and 95\% quantiles are displayed as well. Moreover, the conditional expectations of downward and upward deviations are included, indicating the magnitude of expected deviations at any given time if deviations occur.
Overall, the variance in prosumption is significantly higher during the day compared to nighttime, reflected in broader probabilistic \glspl{ds} during the day across all three cases.

In the following, the results for the three cases are analyzed separately.

\subsubsection{Case 1}

The objective is to minimize residential electricity costs by maximizing self-sufficiency.
Since the \gls{bess} is not used to compensate for uncertainties,
the minimum and maximum battery states coincide, resulting in a deterministic \gls{bess} schedule. This behavior mirrors the underlying framework of \cite{su2021optimal}, where all uncertainties are expected to be compensated by the grid.


Consequently, the probabilities of upward and downward grid deviations always sum up to one. From 09:00 to 22:00, these probabilities stay around 50\%. However, at the start of the schedule (06:00), upward deviations occur with a probability of only 10\% (and consequently downward deviations with around 90\%). This behavior results directly from the highly asymmetrical \gls{pdf} shown in \autoref{fig:pdf-prosumption-forecast} for 06:00. The \gls{pdf} also explains why the expected upward deviation at 06:00 in \autoref{fig:case1}b is particularly high (3 kW), even though its occurrence is unlikely (10\% probability).

At the start of the schedule, the \gls{bess} is discharged, because for midday, dominating \gls{pv} generation is expected. This is incentivized by the quadratic dependence of the nominal grid power $p_G^2$ in the objective, which helps to prevent discharging peaks later in the day. This strategy allows the \gls{bess} to fully utilize its capacity to accommodate the expected high \gls{pv} generation, leading to a plateau in the nominal grid power, see \autoref{fig:case1}b.

In this study, \gls{bess} degradation processes are neglected for simplicity, leading to an optimal strategy of fully discharging and then fully charging the \gls{bess}. To obtain a more accurate solution, \gls{bess} degradation can be directly embedded into the objective function. For example, this could be achieved by integrating the expected \gls{bess} power into a \gls{bess} degradation model, as shown in \cite{su2021optimal} for a deterministic \gls{bess} schedule.

\subsubsection{Case 2}

The primary objective is to minimize residential costs by maximizing self-sufficiency. However, the non-zero weights $c_3$ and $c_4$ reflect that a secondary objective is to minimize uncertainties in grid exchange.

The primary objective translates to first discharging and then fully charging the \gls{bess}, or more precisely, operating the nominal battery state from 0\% to 100\% capacity. For the nominal battery state to reach 0\% (or 100\%) capacity, no upward (or downward) grid uncertainties can be shifted toward the \gls{bess}. Only after the corresponding \gls{bess} limit has been reached, the \gls{bess} can be used for the secondary objective, i.e., to reduce grid uncertainties. In this case, after the nominal battery state reaches 0\% capacity at 09:00, the \gls{bess} is used to compensate for upward grid uncertainties, i.e., the green area beneath the nominal battery state in \autoref{fig:case2}c. After leaving 100\% capacity at 21:00, the \gls{bess} is used to compensate for downward grid uncertainties, i.e., the green area above the nominal battery state.

As a result, mostly upward grid uncertainties are compensated by the \gls{bess}, i.e., the probabilities of upward grid deviations stay around 30\%, whereas the probabilities of downward grid deviations fluctuate around 50\%. This asymmetric compensation arises from the truncation of \glspl{pdf} according to \autoref{fig:pdf-battery}. The maximum battery state aligns with the nominal battery state for the hours 06:00 to 21:00, which indicates that $\ubar{x} = 0$ during that timeframe. Consequently, because the distance between the minimum battery state and the nominal battery state grows over time starting from 09:00, $\Bar{x} \neq 0$ holds. With that, it follows that the expectation of battery power uncertainty $\mathbb{E}[\Delta P_{L\rightarrow B}] \neq 0$ (see \Cref{eq:exp-battery}) and therefore the expected battery state differs from the nominal battery state. To summarize, the asymmetric uncertainty compensation ($\ubar{x} =0, \Bar{x} \neq  0$) results in the difference between the nominal and expected battery state.

\subsubsection{Case 3}

The objective weights for Case 3 portray a tradeoff between minimizing residential costs by increasing self-sufficiency and reducing downward grid uncertainties between 12:00 and 14:00. 

In the resulting probabilistic \gls{ds}, the probabilities of upward grid deviations fluctuate around 40\%. The probabilities of downward deviations during the critical timeframe drop to 20\%, and the 5\% quantiles are closer to the nominal grid power compared to the other cases. However, the nominal battery state cannot be operated to the maximum battery limit as some space is reserved for potential low prosumption realization (e.g., higher \gls{pv} generation), essentially providing flexibility for that critical timeframe.

\subsection{Discussion}
\label{sec:discussion}

In this section, the results of the three cases are discussed together, highlighting key insights of the proposed model. At first, a comparison between cases 1 and 2 is drawn. Afterward, cases 2 and 3 are compared to each other.

\subsubsection{Comparison between Case 1 and Case 2}

In both cases, the nominal grid powers (and consequently also the nominal battery states) align. This indicates that both cases commit to the same nominal grid power schedule which has the same costs for the residents. Additionally, the downward deviations from the nominal grid powers are relatively similar, i.e., the probabilities of downward deviations stay around 50\% for the majority of the day.

However, in Case 2, the probabilities of upward grid deviations are mainly reduced to 30\%, as compared to 50\% in Case 1. This indicates that with 20\% probability, neither downward nor upward deviations occur. This is a key advantage of the proposed model. 
By incorporating a mixed random variable to represent grid power uncertainties (see \autoref{fig:pdf-grid}), the occurrence of discrete events with non-zero probability is enabled:
In 20\% of the cases, the nominal grid power can be adhered to. Those are the cases where the \gls{bess} scheduling accounts for the full compensation of differences between expected and actual prosumption. Since both nominal grid powers align, this is achieved without having to make any compromises regarding the minimization of electricity costs for the homeowner. As a result, an additional degree of freedom in residential \gls{bess} scheduling is identified: The tradeoff lies between reducing grid uncertainties and managing an unknown \gls{bess} state at the end of the schedule---a property further discussed in \Cref{sec:limitations}.

Overall, this comparison demonstrates the model's ability to quantify and reduce grid uncertainties by exploiting the \gls{bess}'s capacity, thereby enabling flexibility without additional electricity costs for the residents.

\subsubsection{Comparison between Case 2 and Case 3}

Case 3 extends the strategy used in Case 2 by actively providing flexibility during a critical period. However, this results in a tradeoff, as the nominal grid power in Case 3 reaches a lower plateau (-3 kW) compared to Case 2 (-2.5 kW), reducing the household's self-sufficiency.

In both cases, the \gls{bess} is fully utilized to reduce grid uncertainties, i.e., the minimum and maximum battery states reach the battery limits. The key difference lies in which uncertainties are compensated at which time. The areas above/below the nominal battery state can be seen as some sort of resource for decreasing downward/upward grid uncertainties, respectively.
Thus, by forcing the \gls{bess} to account for downward grid uncertainties during the critical timeframe, the nominal battery state cannot reach the \gls{bess}'s maximum capacity and thus cannot compensate for as many upward grid uncertainties as in Case 2.

Ultimately, operators must decide, in agreement with residents, how much interference (if any at all) with the optimal residential \gls{bess} schedule is acceptable in order to support the power grid appropriately. Such interference could also be compensated financially, by paying residents for flexibility provision.

\subsection{Limitations}
\label{sec:limitations}

A limitation of the proposed methodology is the absence of a closed-form \gls{pdf} to describe the probabilistic battery state, as discussed in \Cref{sec:battery-model}. While it is not possible to derive a function that fully captures the energy state due to the complexity of summing dependent random variables, key features of the distribution can still be derived, i.e., the minimum, maximum, nominal, and expected battery state. These states provide valuable insight, even though a full analytical representation of the \gls{pdf} remains unavailable for the moment.

Furthermore, it can be noted that the minimum and maximum battery states typically reach the battery limits by the end of the optimization horizon. While grid uncertainties are penalized, uncertainties within the \gls{bess} are not, encouraging full use of the \gls{bess}'s capacity. This however does not guarantee a good starting point for upcoming scheduling tasks, especially if the prosumption is systematically over- or underestimated. Nevertheless, due to the characteristics of the presented model, the expected battery state at the end of the horizon will always be somewhere reasonable around half of the \gls{bess}'s maximum capacity. If this outcome is still undesirable, it can be mitigated by adjusting the \gls{bess} constraints or introducing penalty terms to limit \gls{bess} exploitation.

Additionally, the model allows the \gls{bess} to partially compensate for uncertainties in prosumption. The larger the \gls{bess}'s maximal capacity, the more uncertainties it can absorb. Furthermore, the tighter the prosumption \glspl{pdf}, the easier it is for the \gls{bess} to absorb uncertainties. However, it has been observed that during nighttime periods, the \glspl{pdf} can become too tight. When the weights $\omega_i$ in \Cref{eq:cdf-formula} become excessively small or large, overflow errors may occur during the computation of expected values, causing the optimization process to terminate early.
This issue, however, can be addressed in several ways. One straightforward solution is to broaden the \glspl{pdf} for nighttime hours to make the problem feasible. This is justifiable since in general, \gls{ds} deviations during nighttime are less crucial for the grid anyway and the resulting probabilistic \gls{ds} would simply be a bit more conservative during nighttime. Alternatively, problematic \glspl{pdf} (i.e., \glspl{pdf} with an even steeper increase than the ones presented in \autoref{fig:pdf-prosumption-forecast}) could be approximated using uniform distributions to avoid overflow errors.

Finally, it is important to recognize that forecasting prosumption essentially involves predicting human behavior. It is therefore a highly individual task and forecasting templates, as they are investigated for automating wind and \gls{pv} forecasting \cite{meisenbacher2023autopv}, are assumed to not be easily applicable to the presented methodology. 
Moreover, when looking at the probabilistic \glspl{ds}, it is important to note that exceptions outside of the 95\% quantile can and also will occur. Nevertheless, if those exceptions accumulate, good forecasting models should be able to catch up and adapt.

\section{Conclusion}
\label{sec:conclusion}

We develop a stochastic continuous nonlinear optimization framework for residential day-ahead \acrfull{bess} scheduling. What distinguishes this work from others is the asymmetric allocation of quantified prosumption uncertainties between a residential \gls{bess} and the power grid. This is achieved by introducing mixed random variables to model both \gls{bess} and grid power uncertainties. 
With that, we place great emphasis not only on the \gls{bess} operation but also on the power exchange with the grid.
The model is implemented and tested across several cases, demonstrating that residential \gls{bess} scheduling can effectively reduce grid uncertainties without compromising electricity costs for homeowners. Moreover, during critical periods, the \gls{bess} can be used proactively to mitigate grid uncertainties, contributing to grid stability.
In summary, the proposed model empowers prosumers to actively support future grid stability by managing uncertainties and offering flexibility. This presents a potential step toward integrating residential \glspl{bess} as valuable assets in modern power grids.

Future work will focus on extending the model to intra-day scheduling to enable a better comparison to other scheduling approaches. Additionally, the model's performance in real-world research environments and the interaction between multiple buildings using this approach will be explored to assess its broader impact on grid stability.

\section*{Acknowledgments}
The authors thank the Helmholtz Association for the support under the "Energy System Design" program and the German Research Foundation as part of the Research Training Group 2153 "Energy Status Data: Informatics Methods for its Collection, Analysis and Exploitation".

\section*{Declaration}
For this work, the authors used the free-of-charge ChatGPT version during the writing process to improve clarity and fluency. After the usage, the authors reviewed and edited the content as needed and take full responsibility for the content of the published article. 

\bibliographystyle{elsarticle-num} 
\bibliography{main}

\end{document}